\documentclass{article}

\usepackage{euscript}
\usepackage{amsthm}
\usepackage{amssymb}
\usepackage{amsxtra}

\newtheorem{thm}{Theorem}
\newtheorem{lemma}{Lemma}

\newtheorem{conj}{Conjecture}

\title{A note on Primes in Short Intervals}
\author{Tsz Ho Chan}

\begin{document}
\maketitle
\begin{abstract}
Instead of a strong quantitative form of the Hardy-Littlewood prime $k$-tuple conjecture, one can assume an average form of it and still obtains the same distribution result on $\psi(x+h) - \psi(x)$ by Montgomery and Soundararajan [\ref{MS}].
\end{abstract}
\section{Introduction}
Let $\Lambda(n)$ be von Mangoldt lambda function, $\mu(n)$ be the M\"{o}bius function and $\phi(n)$ be the Euler's phi function. Let $e(\theta) = e^{2\pi i \theta}$, $\epsilon > 0$ and $\psi(x) = \sum_{n \leq x} \Lambda(n)$.
In [\ref{MS}], Montgomery and Soundararajan studied the distribution of primes in short intervals
\begin{equation}
\label{moment}
M_K(N;H) := \sum_{n=1}^{N} (\psi(n+H) - \psi(n) - H)^K
\end{equation}
under a strong quantitative form of the Hardy-Littlewood prime $k$-tuple conjecture:
\begin{conj}
\label{HL1}
\begin{equation}
\label{1}
\sum_{n \leq x} \prod_{i=1}^{k} \Lambda(n+d_i) = {\mathfrak S}(\EuScript{D}) x
+ E_k(x;\EuScript{D})
\end{equation}
holds with
$$E_k(x;\EuScript{D}) \ll_{\epsilon, K} N^{1/2 + \epsilon}$$
uniformly for $1 \leq k \leq K$, $0 \leq x \leq N$, and distinct $d_i$ satisfying $1 \leq d_i \leq H$. Here $\EuScript{D} = \{d_1, d_2, ..., d_k\}$, and
$$\mathfrak{S}(\EuScript{D}) = \mathop{\sum_{q_1, ..., q_k}}_{1 \leq q_i < \infty}
\Bigl(\prod_{i=1}^{k} \frac{\mu(q_i)}{\phi(q_i)} \Bigr) \mathop{\mathop{\mathop{\sum_{a_1, ..., a_k}}_{1 \leq a_i \leq q_i}}_{(a_i, q_i) =1}}_{\sum a_i/q_i \in \mathbb{Z}} e\Bigl(\sum_{i=1}^{k} \frac{a_i d_i}{q_i}\Bigr)$$
is the singular series as in equation (2) of [\ref{MS}].
\end{conj}
They proved
\begin{thm}
\label{thm1}
Under Conjecture \ref{HL1},
\begin{equation*}
\begin{split}
M_K(N;H) =& \mu_k H^{K/2} \int_{1}^{N} (\log x/H + B)^{K/2} dx \\
&+ O\Bigl(N(\log N)^{K/2} H^{K/2} \Bigl(\frac{H}{\log N}\Bigr)^{-1/(8K)} + H^K N^{1/2+\epsilon} \Bigr)
\end{split}
\end{equation*}
uniformly for $\log N \leq H \leq N^{1/K}$, where $\mu_k = 1 \cdot 3 \cdot \cdot \cdot (k-1)$ if $k$ is even, and $\mu_k = 0$ if $k$ is odd; $B = 1 - C_0 - \log 2\pi$ and $C_0$ denotes Euler's constant.
\end{thm}
The first author of [\ref{MS}] suggested to the present author that Theorem \ref{thm1} is probably true under an average form of the Hardy-Littlewood prime $k$-tuple conjecture, namely:
\begin{conj}
\label{HL}
For $x \geq H$,
\begin{equation*}
\mathop{\mathop{\sum_{d_1, ..., d_k}}_{1 \leq d_i \leq H}}_{d_i \textmd{ distinct}} {E_k(x,\EuScript{D})}^2 \ll_k x^{1+\epsilon} H^k.
\end{equation*}
\end{conj}
Our goal in this paper is to prove Theorem \ref{thm1} under Conjecture \ref{HL}.

\bigskip

{\bf Acknowledgement} The author would like to thank Prof. Hugh Montgomery for suggesting the problem and the American Institute of Mathematics for provision both financially and spatially.
\section{Proof of Theorem \ref{thm1}}
\begin{lemma}
\label{lemma1}
With $\Lambda_0(n) = \Lambda(n) - 1$ and $|A|$ denoting the size of a set $A$,
$$\sum_{n \leq x} \prod_{i=1}^{k} \Lambda_0 (n+d_i) = {\mathfrak S}_0 (\EuScript{D}) x + O(\sum_{\EuScript{J} \subset \{1,2,...,k\}} |E_{|\EuScript{J}|} (x; \EuScript{D}_{\EuScript{J}})|)$$
where
$$\mathfrak{S}_0(\EuScript{D}) = \mathop{\sum_{q_1, ..., q_k}}_{1 < q_i < \infty}
\Bigl(\prod_{i=1}^{k} \frac{\mu(q_i)}{\phi(q_i)} \Bigr) \mathop{\mathop{\mathop{\sum_{a_1, ..., a_k}}_{1 \leq a_i \leq q_i}}_{(a_i, q_i) =1}}_{\sum a_i/q_i \in \mathbb{Z}} e\Bigl(\sum_{i=1}^{k} \frac{a_i d_i}{q_i}\Bigr)$$
and $\EuScript{D}_{\EuScript{J}} = \{d_j\}_{j \in \EuScript{J}}$.
\end{lemma}

Proof: The left hand side
\begin{equation*}
\begin{split}
=& \sum_{n \leq x} \prod_{i=1}^{k} (\Lambda(n+d_i) - 1) = \sum_{n \leq x} \sum_{\EuScript{J} \subset \{1,2,...,k\}} (-1)^{k - |\EuScript{J}|} \prod_{i \in \EuScript{J}} \Lambda(n+d_i) \\
=& \sum_{\EuScript{J} \subset \{1,2,...,k\}} (-1)^{k - |\EuScript{J}|} \sum_{n \leq x} \prod_{i \in \EuScript{J}} \Lambda(n+d_i) \\
=& \sum_{\EuScript{J} \subset \{1,2,...,k\}} (-1)^{k - |\EuScript{J}|} ({\mathfrak S}(\EuScript{D}_{\EuScript{J}}) x + E_{|\EuScript{J}|} (x; \EuScript{D}_{\EuScript{J}})) \\
=& {\mathfrak S}_0 (\EuScript{D}) x + O(\sum_{\EuScript{J} \subset \{1,2,...,k\}} |E_{|\EuScript{J}|} (x; \EuScript{D}_{\EuScript{J}})|)
\end{split}
\end{equation*}
by (\ref{1}) and the identity ${\mathfrak S}_0(\EuScript{D}) = \sum_{\EuScript{I} \subset \EuScript{D}} (-1)^{|\EuScript{I}|} {\mathfrak S}(\EuScript{I})$ (see equation (5) of [\ref{MS}]).

\bigskip

Proof of Theorem \ref{thm1} under Conjecture \ref{HL}: Following [\ref{MS}], we expand (\ref{moment}) and have
\begin{equation}
\label{expand}
\begin{split}
M_K(N;H) =& \sum_{k=1}^{K} \frac{1}{k!} \mathop{\sum_{1 \leq M_1, ..., M_k}}_{\sum M_i = K} \binom{K}{M_1 \cdot \cdot \cdot M_k} \\
&\times \mathop{\sum_{m_1,...,m_k}}_{0 \leq m_i < M_i} \prod_{i=1}^{k} (-1)^{M_i - 1 - m_i} \binom{M_i - 1}{m_i} L_k({\bf m}),
\end{split}
\end{equation}
where
\begin{equation}
\label{Ldefine}
L_k({\bf m}) := \mathop{\mathop{\sum_{d_1, ..., d_k}}_{1 \leq d_i \leq H}}_{d_i \textmd{ distinct}} \sum_{n=1}^{N} \prod_{i=1}^{k} \Lambda_{m_i}(n+d_i)
\end{equation}
and $\Lambda_m(n) := \Lambda(n)^m \Lambda_0(n)$. To estimate $L_k({\bf m})$, one needs to distinguish between those $i$ for which $m_i = 0$ and those for which $m_i > 0$. Following [\ref{MS}], we set $\EuScript{K} = \{1, ..., k\}$ and introduce
$$\EuScript{H} = \{i \in \EuScript{K} : m_i \geq 1\}, \; \EuScript{I} = \{i \in \EuScript{K} :  m_i = 0\}, \; \EuScript{J} \subset \EuScript{K}.$$
Then
\begin{equation*}
\begin{split}
& \sum_{n \leq x} \prod_{i \in \EuScript{I}} \Lambda_0(n+d_i) \prod_{i \in \EuScript{H}} \Lambda(n+d_i) \\
=& \sum_{n \leq x} \mathop{\sum_{\EuScript{J}}}_{\EuScript{I} \subset \EuScript{J} \subset \EuScript{K}} \prod_{i \in \EuScript{J}} \Lambda_0(n+d_i) = \mathop{\sum_{\EuScript{J}}}_{\EuScript{I} \subset \EuScript{J} \subset \EuScript{K}} \sum_{n \leq x} \prod_{i \in \EuScript{J}} \Lambda_0(n+d_i) \\
=& \mathop{\sum_{\EuScript{J}}}_{\EuScript{I} \subset \EuScript{J} \subset \EuScript{K}} \Bigl[{\mathfrak S}_0 (\EuScript{D}_{\EuScript{J}}) x + O(\sum_{\EuScript{J}' \subset \EuScript{J}} |E_{|\EuScript{J}'|} (x; \EuScript{D}_{\EuScript{J}'})|) \Bigr] \\
=& x \mathop{\sum_{\EuScript{J}}}_{\EuScript{I} \subset \EuScript{J} \subset \EuScript{K}} {\mathfrak S}_0 (\EuScript{D}_{\EuScript{J}}) + O_k \Bigl( \sum_{\EuScript{K}' \subset \EuScript{K}} |E_{|\EuScript{K}'|} (x; \EuScript{D}_{\EuScript{K}'})| \Bigr)
\end{split}
\end{equation*}
by Lemma \ref{lemma1}. We write the above as $f(x) = cx + E_{x,\EuScript{K}}$. In general,
\begin{equation*}
\begin{split}
& \int_{1^-}^{X} g(x) df(x) \\
=& g(X) f(X) - \int_{1^-}^{X} f(x) g'(x) dx \\
=& g(X)[cX + E_{X,\EuScript{K}}] - \int_{1^-}^{X} (cx + E_{x,\EuScript{K}}) g'(x) dx \\
=& c \int_{1^-}^{X} g(x) dx + O\Bigl(E_{X,\EuScript{K}}|g(X)| + \int_{1^-}^{X} E_{x,\EuScript{K}} |g'(x)| dx  \Bigr).
\end{split}
\end{equation*}
Thus, by integration by parts,
\begin{equation*}
\begin{split}
\sum_{n=1}^{N} & \Bigl(\prod_{i \in \EuScript{I}} \Lambda_0(n+d_i)\Bigr) \Bigl(\prod_{i \in \EuScript{H}} \Lambda(n+d_i) (\log (n+d_i))^{m_i -1} (\log (n+d_i) - 1)\Bigr) \\
=& \int_{1^-}^{N} \prod_{i \in \EuScript{H}} (\log (x+d_i))^{m_i - 1} (\log (x+d_i) - 1) df(x) \\
=&c \int_{1}^{N} \prod_{i \in \EuScript{H}} (\log{(x+d_i)})^{m_i - 1} (\log{(x+d_i)} - 1) dx + O\Bigl(E_{N,\EuScript{K}} \log^K N \\
&+ \int_{1^-}^{N} E_{x,\EuScript{K}} \frac{\log^K (x+H)}{x} dx \Bigr).
\end{split}
\end{equation*}
This is the analogue of equation (65) in [\ref{MS}]. Now, note that
$$\Lambda_m(n) = \Lambda(n) (\log n)^{m-1} (\log n -1)$$
when $n$ is prime. We have, by following the argument in [\ref{MS}],
\begin{equation}
\label{middle}
\begin{split}
& \sum_{n=1}^{N} \prod_{i=1}^{k} \Lambda_{m_i}(n+d_i) \\
=& \Bigl(\mathop{\sum_{\EuScript{J}}}_{\EuScript{I} \subset \EuScript{J} \subset \EuScript{K}} {\mathfrak S}_0 (\EuScript{D}_{\EuScript{J}}) \Bigr) (I_{\bf m}(N) + O(H(\log N)^{K-k})) + O(N^{1/2+\epsilon}) \\
&+ O\Bigl(E_{N,\EuScript{K}} \log^K N + \int_{1^-}^{N} E_{x,\EuScript{K}} \frac{\log^K (x+H)}{x} dx \Bigr).
\end{split}
\end{equation}
where
$$I_{\bf m}(N) := \int_{1}^{N} \prod_{i \in \EuScript{H}} \bigl((\log x)^{m_i - 1} (\log x - 1)\bigr) dx.$$
Putting (\ref{middle}) into (\ref{Ldefine}),
\begin{equation}
\label{L}
\begin{split}
L_k({\bf m}) =& I_{\bf m}(N) \mathop{\sum_{\EuScript{J}}}_{\EuScript{I} \subset \EuScript{J} \subset \EuScript{K}} \mathop{\mathop{\sum_{d_1, ..., d_k}}_{1 \leq d_i \leq H}}_{d_i \textmd{ distinct}} {\mathfrak S}_0 (\EuScript{D}_{\EuScript{J}}) + O(H^k N^{1/2+\epsilon}) \\
&+ O\Bigl(\mathop{\mathop{\sum_{d_1, ..., d_k}}_{1 \leq d_i \leq H}}_{d_i \textmd{ distinct}} E_{N,\EuScript{K}} \log^K N + \int_{1^-}^{N} \mathop{\mathop{\sum_{d_1, ..., d_k}}_{1 \leq d_i \leq H}}_{d_i \textmd{ distinct}} E_{x,\EuScript{K}} \frac{\log^K x}{x} dx \Bigr).
\end{split}
\end{equation}
Now, we use Conjecture \ref{HL}. By Cauchy-Schwarz inequality,
\begin{equation*}
\label{average}
\mathop{\mathop{\sum_{d_1, ..., d_k}}_{1 \leq d_i \leq H}}_{d_i \textmd{ distinct}} |E_k(x,\EuScript{D})| \ll_k x^{1/2+\epsilon} H^k.
\end{equation*}
In particular,
\begin{equation*}
\begin{split}
\mathop{\mathop{\sum_{d_1, ..., d_k}}_{1 \leq d_i \leq H}}_{d_i \textmd{ distinct}} E_{x,\EuScript{K}} =& \mathop{\mathop{\sum_{d_1, ..., d_k}}_{1 \leq d_i \leq H}}_{d_i \textmd{ distinct}} \sum_{\EuScript{K}' \subset \EuScript{K}} |E_{|\EuScript{K}'|} (x; \EuScript{D}_{\EuScript{K}'})| \\
&\ll_k \sum_{j=0}^{k} \mathop{\mathop{\sum_{d_1, ..., d_j}}_{1 \leq d_i \leq H}}_{d_i \textmd{ distinct}} |E_j(x;\EuScript{D})| \ll_k x^{1/2+\epsilon} H^k
\end{split}
\end{equation*}
Applying this to (\ref{L}), the second error term is $\ll N^{1/2 + \epsilon} H^k$ while the third error term is
\begin{equation*}
\begin{split}
=& \int_{H}^{N} \mathop{\mathop{\sum_{d_1, ..., d_k}}_{1 \leq d_i \leq H}}_{d_i \textmd{ distinct}} E_{x,\EuScript{K}} \frac{\log^K (x+H)}{x} dx + \int_{1^-}^{H} \mathop{\mathop{\sum_{d_1, ..., d_k}}_{1 \leq d_i \leq H}}_{d_i \textmd{ distinct}} E_{x,\EuScript{K}} \frac{\log^K (x+H)}{x} dx \\
\ll_k& N^{1/2+\epsilon} H^k + \int_{1}^{H} H^{k+1} \frac{\log^K (x+H)}{x} dx \ll_K N^{1/2+\epsilon} H^k
\end{split}
\end{equation*}
as $H \leq N^{1/2}$. Hence, (\ref{L}) has an error $O(H^k N^{1/2+\epsilon})$ and the rest of the proof in [\ref{MS}] follows. Therefore, we have Theorem \ref{thm1} under Conjecture \ref{HL}.


Tsz Ho Chan\\
American Institute of Mathematics\\
360 Portage Avenue\\
Palo Alto, CA 94306\\
U.S.A.\\
thchan@aimath.org

\end{document}